\documentclass[12pt]{amsart}

\usepackage{amsthm}

\usepackage{mathrsfs,mathtools,latexsym,eucal}
\usepackage{amscd,amsfonts,amssymb,amsmath,amsthm}
\usepackage{mathrsfs}
\usepackage{graphicx,graphics,color}
\usepackage[utf8]{inputenc}
\usepackage[active]{srcltx}
\usepackage[english]{babel}
\usepackage[pagewise,displaymath,mathlines]{lineno}
\usepackage{listings}
\usepackage{color}
\usepackage{subfigure} 
\usepackage{url}
\usepackage{float}
\usepackage{enumerate}

\usepackage[active]{srcltx}

\usepackage{epstopdf}

\usepackage{pgf,tikz}
\usetikzlibrary{arrows}

\definecolor{dkgreen}{rgb}{0,0.6,0}
\definecolor{gray}{rgb}{0.5,0.5,0.5}
\definecolor{mauve}{rgb}{0.58,0,0.82}

\newtheorem*{theorem*}{Theorem}
\newtheorem{theorem}{Theorem}[section]
\newtheorem{lemma}[theorem]{Lemma}

\newtheorem{corollary}[theorem]{Corollary}

\newtheorem{remark}[theorem]{Remark}

\newtheorem{definition}[theorem]{Definition}

\newcommand{\C}{\mathbb{C}}
\newcommand{\R}{\mathbb{R}}
\newcommand{\N}{\mathbb{N}}

\renewcommand{\leq}{\leqslant}
\renewcommand{\geq}{\geqslant}

\numberwithin{equation}{section}

\newcommand{\Bo}{\mathcal{B}}

\newcommand{\Un}{\mathcal{U}}

\newcommand{\spa}{\mathrm{span}}

\newcommand{\suppv}{\mathrm{suppv}}

\newcommand{\gsuppv}{\mathrm{gsuppv}}

\newcommand{\B}{\mathsf{B}}

\newcommand{\E}{\mathsf{S}}

\newcommand{\argmax}{\mathrm{arg}\max}
\newcommand{\argmin}{\mathrm{arg}\min}

\begin{document}

\title[Exact solutions to $\displaystyle{\max_{\|x\|=1} \sum_{i=1}^\infty\|T_i(x)\|^2}$ with applications to Physics, Bioengineering and Statistics]{Exact solutions to $\displaystyle{\max_{\|x\|=1} \sum_{i=1}^\infty\|T_i(x)\|^2}$ with applications to Physics,  Bioengineering and Statistics}

\author{Francisco Javier Garcia-Pacheco$^0$}
\address{Department of Mathematics, College of Engineering, University of Cadiz, Puerto Real 11510, Spain (EU)}
\email[Corresponding author]{{\tt garcia.pacheco@uca.es}}

\author{Clemente Cobos-Sanchez}
\address{Department of Electronics, College of Engineering, University of Cadiz, Puerto Real 11510, Spain (EU)}
\email{{\tt clemente.cobos@uca.es}}

\author{Soledad Moreno-Pulido$^0$}
\address{Department of Mathematics, College of Engineering, University of Cadiz, Puerto Real 11510, Spain (EU)}
\email{{\tt soledad.moreno@uca.es}}

\author{Alberto Sanchez-Alzola}
\address{Department of Statistics and Operation Research, College of Engineering, University of Cadiz, Puerto Real 11510, Spain (EU)}
\email{{\tt alberto.sanchez@uca.es}}

\date{}

\footnotetext{The authors have been supported by the Research Grant PGC-101514-B-100 awarded by the Spanish Ministry of Science, Innovation and Universities.}
\thanks{The authors would like to express their deepest gratitude towards the reviewers for their valuable suggestions and comments that helped improve the paper considerably. All authors contributed equally to the work}

\begin{abstract}
The supporting vectors of a matrix $A$ are the solutions of $\displaystyle{\max_{\|x\|_2=1} \|Ax\|_2^2}$. The generalized supporting vectors of matrices $A_1,\dots,A_k$ are the solutions of the problem $\displaystyle{\max_{\|x\|_2=1} \|A_1x\|_2^2+\cdots+\|A_kx\|_2^2}$. Notice that the previous optimization problem is also a boundary element problem since the maximum is attained on the unit sphere. Many problems in Physics, Statistics and Engineering can be modeled by using generalized supporting vectors. In this manuscript we first raise the generalized supporting vectors to the infinite dimensional case by solving the optimization problem $\displaystyle{\max_{\|x\|=1} \sum_{i=1}^\infty\|T_i(x)\|^2}$ where $(T_i)_{i\in\N}$ is a sequence of bounded linear operators between Hilbert spaces $H$ and $K$ of any dimension. Observe that the previous optimization problem generalizes the first two. Then a unified MATLAB code is presented for computing generalized supporting vectors of a finite number of matrices. Some particular cases are considered and three novel examples are provided to which our technique applies: optimized observable magnitudes by a pure state in a quantum mechanical system, a TMS optimized coil and an optimal location problem using statistics multivariate analysis. These three examples show the wide applicability of our theoretical and computational model.
\end{abstract}

\maketitle

\section{Introduction}

Many problems in different disciplines like Physics, Statistics, Economics or Engineering can be modeled by using matrices and their norms (see for instance \cite{H,Y}). These real-life problems usually look like:
\begin{equation}\label{original}
\left\{\begin{array}{l}
    \max \|T_i(x)\|\;\;\; i\in \N\\
    \min\|x\|
    \end{array}\right.
\end{equation}
where $T_i:H\to K$ is a continuous or bounded linear operator between real or complex Hilbert spaces $H$ and $K$ for every $i\in\N$. In most situations, $H$ and $K$ are both finite dimensional real Hilbert spaces and we have a finite number of operators (matrices), that is, there exists $m\in \N$ such that $T_i=0$ if $i\geq m$. However, there still are real-life situations where $H=K$ is an infinite dimensional complex Hilbert space (see Subsection \ref{quantum}).

The first aspect to take into consideration about Problem \eqref{original} is the fact that it has no solution unless $T_i=0$ for every $i\in\N$. Indeed, it is a multiobjective optimization problem with the only constraint that $x\in H$, so if $x_0\in H$ solves \eqref{original}, then in particular $x_0$ minimizes $\|x\|$ so $x_0=0$. Since $x_0$ maximizes $\|T_i(x)\|$ for every $i\in\N$, this forces that $T_i=0$ for every $i\in \N$, which in general does not occur. Furthermore, it is also useful to bear in mind that $\max \|T(x)\|$, for $T$ a nonzero bounded operator from $H$ to $K$, is an optimization problem with no solution at all, since if $x_0\in H$ is a solution of $\max\|T(x)\|$, then it is clear that $x_0\notin \ker(T)$ and thus $kx_0$, $k\in \N$, gives greater values of $\|T(\bullet)\|$ than $x_0$. In fact, $\|T(kx_0)\|=k\|T(x_0)\|\to \infty$ as $k\to \infty$.

The above reasons imply that Problem \eqref{original} must be appropriately reformulated in such a way that the following conditions must be met:
\begin{enumerate}
    \item The reformulation must be consistent with the real-life problem.
    \item The set of solutions of the reformulation must be nonempty.
    \item In the finite dimensional case with a finite number of matrices, finding the solution must be computationally affordable.
\end{enumerate}

In the following section, we will recall some basic concepts of Optimization Theory that will help us reformulate Problem \eqref{original}. However, it seems clear that in order to meet the above conditions we will have to introduce a constraint and gather all the objective functions to maximize into a single objective function to maximize.

\section{Preliminaries}

A multiobjective optimization problem has the form $$P:=\left\{\begin{array}{l} \min f_j(x),\; 1\leq j \leq l \\ g_i(x) \leq b_i,\; 1\leq i\leq k\end{array}\right.$$ where $f_1,\dots,f_l,g_1,\dots,g_k: X\to \mathbb{R}$ are functions defined on a nonempty set $X$. Two special sets are associated to $P$, the feasible solutions of $P$ $$\mathrm{fea}(P):=\{x\in X: g_i(x) \leq b_i\; \forall\ 1\leq i\leq k\}$$ and the set of optimal solutions of $P$ $$\mathrm{sol}(P):=\{x\in \mathrm{fea}(P): f_j(x)\leq f_j(y)\; \forall y\in \mathrm{fea}(P)\;\forall\ 1\leq j\leq l\}.$$

Any multiobjective optimization problem can be rewritten as the intersection of optimization problems, that is, if $$P_j:=\left\{\begin{array}{l} \min f_j(x) \\ g_i(x) \leq b_i,\; 1\leq i\leq k\end{array}\right.$$ for $1\leq j\leq l$, then $P=P_1\wedge \cdots \wedge P_l$, $\mathrm{fea}(P)=\mathrm{fea}(P_j)$ for all $1\leq j\leq l$, and $\mathrm{sol}(P)=\mathrm{sol}(P_1)\cap \cdots \cap \mathrm{sol}(P_l)$. In case that $\mathrm{sol}(P)=\mathrm{sol}(P_1)\cap \cdots \cap \mathrm{sol}(P_l)=\varnothing,$ the multiobjective optimization problem must be reformulated. A typical reformulation is 
\begin{equation}\label{refor1}
\left\{\begin{array}{l} \min f_j(x),\; 1\leq j \leq l \\ g_i(x) \leq b_i,\; 1\leq i\leq k\end{array}\right.\stackrel{\text{reform}}{\longrightarrow}  \left\{\begin{array}{l} \min f_{j_0}(x)\\ f_j(x)\leq c_j,\; 1\leq j \leq l, \; j\neq j_0\\ g_i(x) \leq b_i,\; 1\leq i\leq k\end{array}\right.
\end{equation}
Another typical reformulation is
\begin{equation}\label{refor2}
\left\{\begin{array}{l} \min f_j(x),\; 1\leq j \leq l \\ g_i(x) \leq b_i,\; 1\leq i\leq k\end{array}\right.\stackrel{\text{reform}}{\longrightarrow}  \left\{\begin{array}{l} \min h\left(f_1(x),\dots,f_l(x)\right)\\  g_i(x) \leq b_i,\; 1\leq i\leq k\end{array}\right.
\end{equation}
where $h:\R^l\to\R$ is a function conveniently chosen (usually an increasing function on each coordinate).

On the other hand, observe that if $\phi:Y\to X$ is a bijection, then it is easy to check that $\mathrm{fea}(P)=\phi(\mathrm{fea}(Q))$ and $\mathrm{sol}(P)=\phi(\mathrm{sol}(Q))$ where
\begin{equation}\label{cv1}
Q:=\left\{\begin{array}{l} \min (f_j\circ\phi)(y),\; 1\leq j \leq l \\ (g_i\circ \phi)(y) \leq b_i,\; 1\leq i\leq k\end{array}\right.
\end{equation}

Also note that $\mathrm{fea}(P)=\mathrm{fea}(R)$ and $\mathrm{sol}(P)=\mathrm{sol}(R)$ where
\begin{equation}\label{cv2}
R:=\left\{\begin{array}{l} \min (\phi_j\circ f_j)(x),\; 1\leq j \leq l \\ (\chi_i\circ g_i)(x) \leq \chi_i(b_i),\; 1\leq i\leq k\end{array}\right.
\end{equation}
and $\phi_j,\chi_i: \mathbb{R}\to \mathbb{R}$ are strictly increasing for $1\leq j\leq l$ and $1\leq i\leq k$.

\section{Reformulation of Problem \eqref{original}}\label{orgsec}

Let us get back to Problem \eqref{original}. Like we have previously mentioned in the first section, it is convenient to gather all the objective functions to maximize into a single objective function to maximize. This means that we have to come up with a bounded linear operator $T$ that involves in some sense the $T_i$’s. Note that we are working in a Hilbert space setting, and since $$\ell_2(K):=\left\{(y_i)_{i\in\N}\in K^{\N}: \sum_{i=1}^\infty \|y_i\|^2<\infty\right\}$$ endowed with $\|(y_i)_{i\in\N}\|_2:=\left(\sum_{i=1}^\infty \|y_i\|^2\right)^{\frac{1}{2}}$ is a Hilbert space, the sequence $(T_i(x))_{i\in\N}$ should verify that $(T_i(x))_{i\in\N}\in \ell_2(K)$, in other words, $\sum_{i=1}^\infty \|T_i(x)\|^2<\infty$. The idea is not to fall out of the Hilbert space setting. This leads us to the following reformulation
\begin{equation}\label{refor0}
    \left\{\begin{array}{l}
     \max \displaystyle{\sum_{i=1}^\infty \|T_i(x)\|^2}  \\
         \|x\|= 1.
    \end{array}
    \right.
\end{equation}
Here we have applied the first typical reformulation \eqref{refor1} to $\min\|x\|$ and we have obtained the constraint $\|x\|\leq 1$. The homogeneous character of the norm together with the linearity of the operators make irrelevant the choice of any other constant $a>0$ in $\|x\|\leq a$ since we can normalize. For the same reason, the maximum on $\|x\|\leq 1$ is attained on $\|x\|=1$. This is why the constraint $\|x\|=1$ appears in \eqref{refor0}. In the reformulation \eqref{refor0} we have also applied the second typical reformulation \eqref{refor2} to the $\max \|T_i(x)\|$’s and we have obtained $\max \sum_{i=1}^\infty \|T_i(x)\|^2$ taking into consideration that the sequence $(T_i(x))_{i\in\N}\in \ell_2(K)$. Now, if we consider the operator 
\begin{equation}\label{rightoperator}
    \begin{array}{rrcl}
       T:  &  H&\to & \ell_2(K)\\
         & x & \mapsto & T(x):=\left( T_i(x)\right)_{i\in\N},
    \end{array}
\end{equation}
then the reformulation \eqref{refor0} can be rewritten as 
\begin{equation}\label{refor01}
\left\{\begin{array}{l}
    \max \|T(x)\|_2\\
    \|x\|= 1.
 \end{array}   \right.
\end{equation}
However, the operator $T$ given in \eqref{rightoperator} involves the Hilbert space $\ell_2(K)$, which can be hard to handle computationally in the finite dimensional case with a large finite number of matrices, simply because we might have to compute many Cholesky decompositions. In our main theorem (Theorem \ref{gsv}), we solve exactly the reformulation \eqref{refor0} and show that its solutions also solve
\begin{equation}\label{refor02}
\left\{\begin{array}{l}
    \max \|R(x)\|\\
    \|x\|= 1.
 \end{array}   \right.
\end{equation}
where now $R:=\sum_{i=1}^\infty T_i'\circ T_i$, being $T_i':K\to H$ the adjoint operator of $T_i$.

To conclude this section, in order to reinforce the choice of the reformulation \eqref{refor0}, we will prove the following theorem, which states that \eqref{refor0} is equivalent to \begin{equation}\label{refor03}
\left\{\begin{array}{l}
    \min \dfrac{\|x\|^2}{\displaystyle{\sum_{i=1}^\infty \|T_i(x)\|^2}}\\
    \\
    \displaystyle{x\notin \bigcap_{i\in\N}\ker(T_i)}.
 \end{array}   \right.
\end{equation}

Note that the function $$\dfrac{\|x\|^2}{\displaystyle{\sum_{i=1}^\infty \|T_i(x)\|^2}}$$ is homogeneous of degree $0$, which implies that if $x_0$ is a solution of \eqref{refor03}, then $\lambda x_0$ is also a solution for every $\lambda\neq 0$.

\begin{theorem}\label{reforequiv}
Let $H$ and $K$ be Hilbert spaces. Let $(T_i)_{i\in\N}$ be a sequence of bounded linear operators from $H$ to $K$, not all zero, such that $\sum_{i=1}^\infty \|T_i\|^2<\infty$. Then $$\argmax_{\|x\|=1} \sum_{i=1}^\infty \|T_i(x)\|^2 =\E_H\cap \left( \argmin_{x\notin \bigcap_{i\in\N}\ker(T_i)} \dfrac{\|x\|^2}{\displaystyle{\sum_{i=1}^\infty \|T_i(x)\|^2}}\right).$$
\end{theorem}

\begin{proof}
Let $$x_0\in \argmax_{\|x\|=1} \sum_{i=1}^\infty \|T_i(x)\|^2.$$ We will first show that $x_0\notin \bigcap_{i\in\N}\ker(T_i)$. Suppose to the contrary that $x_0 \in \bigcap_{i\in\N}\ker(T_i)$. Then $\sum_{i=1}^\infty \|T_i(x_0)\|^2=0$. Since not all the $T_i$'s are zero, we can find $x_1\in H$ and $i_1\in \N$ such that $T_{i_1}(x_1)\neq 0$, which gives us the following contradiction: $$\sum_{i=1}^\infty \|T_i(x_0)\|^2=0<\sum_{i=1}^\infty \|T_i(x_1)\|^2.$$ Thus $x_0\notin \bigcap_{i\in\N}\ker(T_i)$. Fix an arbitrary $y\in H$ such that $y\notin \bigcap_{i\in\N}\ker(T_i)$. We will show that
\begin{equation}\label{zero}
\dfrac{\|x_0\|^2}{\displaystyle{\sum_{i=1}^\infty \|T_i(x_0)\|^2}}\leq \dfrac{\|y\|^2}{\displaystyle{\sum_{i=1}^\infty \|T_i(y)\|^2}}.
\end{equation}
By hypothesis,
\begin{equation}\label{uno}
\sum_{i=1}^\infty \left\|T_i\left(\frac{y}{\|y\|}\right)\right\|^2\leq \sum_{i=1}^\infty \|T_i(x_0)\|^2.
\end{equation}
By rearranging terms in Equation \eqref{uno} and by bearing in mind that $\|x_0\|=1$, we obtain \eqref{zero}. Conversely, let $$x_0\in\E_H\cap \left(\argmin_{x\notin \bigcap_{i\in\N}\ker(T_i)} \dfrac{\|x\|^2}{\displaystyle{\sum_{i=1}^\infty \|T_i(x)\|^2}}\right).$$ Fix an arbitrary $y\in\E_H$. We will show that
\begin{equation}\label{dos}
\sum_{i=1}^\infty \left\|T_i\left(y\right)\right\|^2\leq \sum_{i=1}^\infty \left\|T_i\left(x_0\right)\right\|^2.
\end{equation}
By hypothesis,
\begin{equation}\label{tres}
\dfrac{\|x_0\|^2}{\displaystyle{\sum_{i=1}^\infty \|T_i(x_0)\|^2}}\leq \dfrac{\|y\|^2}{\displaystyle{\sum_{i=1}^\infty \|T_i(y)\|^2}}.
\end{equation}
By rearranging terms in Equation \eqref{tres} and by bearing in mind that $\|y\|=\|x_0\|=1$, we obtain \eqref{dos}.
\end{proof}

Reformulations of the form given in \eqref{refor03} have been widely considered in the finite dimensional case with a finite number of matrices (see, for instance, \cite{key:wassermann,key:romei}) simply because it is easy to apply an heuristic method to approximate a solution. Our main result (Theorem \ref{gsv}) solves \eqref{refor0} exactly and thus \eqref{refor03} exactly too with no need of applying heuristic methods which many times are not proven to converge to the solution.

\section{Generalized supporting vectors}

Supporting vectors are widely known in the literature of Geometry of Banach Spaces and Operator Theory. They are commonly known as the unit vectors at which an operator attains the maximum of its norm. In the matrix setting, the supporting vectors of a matrix $A$ are the solutions of $$\displaystyle{\max_{\|x\|_2=1} \|Ax\|_2^2}.$$ In \cite{CSGPMPSM,GPNG} supporting vectors are topologically and geometrically studied. In fact, the first relevant papers where supporting vectors of a given fixed operator are studied for the first time are \cite{CSGPMPSM,GPNG,CSGPGRH}. In these papers, the approach is completely different as the ones given in \cite{BP1,BP2,James,Lin}, since in the latter papers the focus is given on the density of the norm-attaining functionals or operators in the dual space whereas in \cite{CSGPMPSM,GPNG,CSGPGRH} the density of those functionals/operators is disregarded and the spotlight is occupied by the supporting vectors instead. The reason why we concentrate on the supporting vectors is because many real-life problems are formulated in terms of maximizing the norm of a given fixed matrix, therefore to solve those optimization problems we need the supporting vectors. As far as we know, no real-life problem is formulated in terms of the density of norm-attaining matrices nor such density is involved at all in solving or setting applied optimization problems.

In addition, in \cite{CSGPMPSM} generalized supporting vectors are defined and studied. Again in the matrix setting, the generalized supporting vectors of a sequence of matrices $(A_i)_{i\in\N}$ are the solutions of $$\displaystyle{\max_{\|x\|_2=1}\sum_{i=1}^\infty \|A_ix\|_2^2}.$$ This optimization problem clearly generalizes the previous one.

Let us go over these concepts with the formalism proper from Functional Analysis.

Let $X$ be a complex Banach space. By $\Bo(X)$ we denote the Banach space of continuous or bounded linear operators on $X$. Notice that in virtue of the Open Mapping Theorem, the group of invertibles of $\Bo(X)$ is $$\Un(\Bo(X)):=\left\{T\in\Bo(X): \ker(T)=\{0\}\text{ and }T(X)=X\right\}.$$ The spectrum of an operator $T\Bo(X)$ is defined as $$\sigma(T):=\left\{\lambda\in\C: T-\lambda I \notin \Un(\Bo(X))\right\}.$$ The spectral decomposition theorem states that the the spectrum is the disjoint union of the point spectrum, the continuous spectrum and the residual spectrum, in other words, $\sigma(T)=\sigma_p(T)\cup \sigma_c(T)\cup \sigma_r(T)$, where the point spectrum is $$\sigma_p(T):=\{\lambda\in\C:\ker(T-\lambda I)\neq \{0\}\},$$ the continuous spectrum is $$\sigma_c(T):=\{\lambda\in\C:\ker(T-\lambda I)= \{0\}\text{ and }T(X)\subsetneq \overline{T(X)}=X\}$$ and the residual spectrum is $$\sigma_p(T):=\{\lambda\in\C:\ker(T-\lambda I)= \{0\} \text{ and }\overline{T(X)}\neq X\}.$$ The elements of the point spectrum are called the eigenvalues of $T$. If $\lambda\in\sigma_p(T)$, then $$V(\lambda):=\ker(T-\lambda I)=\{x\in X: T(x)=\lambda x\}$$ is called the subspace of eigenvectors associated to the eigenvalue $\lambda$.

It is trivial that if $\lambda\in \sigma_p(T)$, then $|\lambda|\leq \|T\|$. Therefore, if $\|T\|\in \sigma_p(T)$, then it is the largest eigenvalue, which is denoted as $\lambda_{\max}(T)$. In this case, $\E_X\cap V(\|T\|)\subseteq \suppv(T)$, where $\E_X$ is the unit sphere of $X$, that is, the set of unit vectors of $X$, and $\suppv(T):=\{x\in\E_X:\|T(x)\|=\|T\|\}$ is the set of supporting vectors of $T$ (see \cite{CSGPMPSM,GPNG}). If $(T_n)_{n\in\N}\subset \Bo(X,Y)$ (the Banach space of continuous or bounded linear operators from $X$ to $Y$), then the set of generalized supporting vectors (see \cite[Definition 3.1]{CSGPMPSM}) of $(T_n)_{n\in\N}$ is defined as $$\gsuppv\left((T_n)_{n\in\N}\right):=\argmax_{\|x\|=1}\sum_{n=1}^\infty\|T_n(x)\|^2.$$ It is not hard to check that $\gsuppv\left((T_n)_{n\in\N}\right)=\suppv\left(S\right)$ where $$\begin{array}{rrcl} S:& X&\to &\ell_2(Y)\\ & x &\mapsto& S(x):= (T_n(x))_{n\in\N}\end{array}$$ and $\ell_2(Y):=\left\{(y_n)_{n\in\N}\in X^\N : \sum_{n=1}^\infty \|y_n\|^2<\infty\right\}$. For the operator $S$ to be well defined is sufficient that $\sum_{n=1}^\infty \|T_n\|^2<\infty$.

In \cite[Theorem 3.3]{CSGPMPSM} it was shown that 
\begin{eqnarray*}
    \max_{\|x\|_2 = 1} \sum_{i=1}^k \| A_i x\|_2^2 &=& \lambda_{\max} \left(\sum_{i=1}^k A_i^T A_i\right) \\
    \arg\max_{\|x\|_2 = 1} \sum_{i=1}^k \| A_i x\|_2^2  &=& V\left(\lambda_{\max}\Big(\sum_{i=1}^k A_i^T A_i\Big)\right) \cap \E_{\ell_2^n}
\end{eqnarray*}

where $A_1,\dots,A_k$ are $m\times n$ real matrices. In this section we intend to generalize this result to the infinite dimensional case.

Recall that if $H$ is a Hilbert space, then the dual map of $H$ is defined as $$\begin{array}{rrcl} J_H: & H&\mapsto & H^*\\ &h&\mapsto &J_H(h):=h^*=(\bullet|h)\end{array}$$ This identification $J_H$ is in fact a surjective linear isometry between $H$ and $H^*$ and in the setting of the Geometry of Banach Spaces is also known as the duality mapping.

Now, if $H$ and $K$ are Hilbert spaces and $T\in \Bo(H,K)$, then the adjoint operator of $T$ is defined as $T’:=(J_H)^{-1}\circ T^*\circ J_K\in \Bo(K,H)$, where $T^*:K^*\to H^*$ is the dual operator of $T$. Among all the properties verified by the adjoint $T’$ we have that $(T(x)|y)=(x|T'(y))$ for all $x\in H$ and all $y\in K$ and $\|T’\|=\|T\|$. An operator $T\in\Bo(H)$ is said to be a selfadjoint operator provided that $T=T'$. It is not hard to check that all the eigenvalues of a selfadjoint operator are real. On the other hand, $T$ is said to be positive provided that $(T(x)|x)\geq 0$ for all $x\in H$. It is trivial that the eigenvalues of a positive operator have to be positive. Finally, $T$ is said to be normal provided that $T=S'\circ S$ for some $S\in\Bo(H,K)$. Normal operators are examples of selfadjoint positive operators.

It is well known that if $T\in\Bo(H)$ is a selfadjoint, positive and compact operator on a Hilbert space $H$, then $\|T\|=\lambda_{\max}(T)$.

Note that if $(T_i)_{i\in\N}\subset \Bo(H,K)$ is a sequence of bounded linear operators such that such that $\sum_{i=1}^\infty\|T_i\|^2<\infty$, then $\sum_{i=1}^\infty T_i'\circ T_i\in \Bo(H)$ is trivially selfadjoint and positive. If, in addition, the $T_i$’s are compact, then $\sum_{i=1}^\infty T_i'\circ T_i$ is also compact. In fact, the compact operators form a closed subalgebra of $\Bo(H)$ commonly denoted as $\mathcal{K}(H)$.

\begin{lemma}\label{gsv0}
Let $H$ and $K$ be Hilbert spaces and $(T_i)_{i\in\N}\subset \Bo(H,K)$ a sequence of bounded linear operators such that $\sum_{i=1}^\infty\|T_i\|^2<\infty$. Then $$\sum_{i=1}^\infty\|T_i(x)\|^2 = x^*\left(\left(\sum_{i=1}^\infty T_i'\circ T_i\right)(x)\right)\leq \left\|\left(\sum_{i=1}^\infty T_i' \circ T_i\right)(x)\right\|\leq \left\|\sum_{i=1}^\infty T_i' \circ T_i\right\|$$ for all $x\in H$.
\end{lemma}

\begin{proof}
Notice that $\sum_{i=1}^\infty T_i'\circ T_i$ is a well defined operator since it is an absolutely convergent series. Indeed, $$\sum_{i=1}^\infty \|T_i'\circ T_i\|\leq \sum_{i=1}^\infty \|T_i'\|\| T_i\|=\sum_{i=1}^\infty \|T_i\|^2<\infty.$$ Also, $\sum_{i=1}^\infty\|T_i(x)\|^2\leq\|x\|^2 \sum_{i=1}^\infty\|T_i\|^2<\infty$ is well defined. Finally, it suffices to observe that

\begin{eqnarray*}
\sum_{i=1}^\infty\|T_i(x)\|^2 &=& \sum_{i=1}^\infty \left(T_i(x)|T_i(x)\right)\\
&=& \sum_{i=1}^\infty \left(T_i'(T_i(x))|x\right)\\
&=& \sum_{i=1}^\infty x^*\left(T_i'(T_i(x))\right)\\
&=& x^*\left(\left(\sum_{i=1}^\infty T_i'\circ T_i\right)(x)\right)\\
&\leq & \|x^*\| \left\|\left(\sum_{i=1}^\infty T_i' \circ T_i\right)(x)\right\|\\
 &=&\left\|\left(\sum_{i=1}^\infty T_i' \circ T_i\right)(x)\right\|\\
 &\leq & \left\|\sum_{i=1}^\infty T_i' \circ T_i\right\|.
\end{eqnarray*}
\end{proof}

The following result is a generalization of \cite[Theorem 3.3]{CSGPMPSM} to the infinite dimensional case.

\begin{theorem}\label{gsv}
Let $H$ and $K$ be Hilbert spaces and $(T_i)_{i\in\N}\subset \Bo(H,K)$ a sequence of compact operators such that $\sum_{i=1}^\infty\|T_i\|^2<\infty$. Then $$\max_{\|x\|=1} \sum_{i=1}^\infty\|T_i(x)\|^2 = \lambda_{\max}\left(\sum_{i=1}^\infty T_i'\circ T_i\right)=\left\|\sum_{i=1}^\infty T_i'\circ T_i\right\|$$ and $$\argmax_{\|x\|=1} \sum_{i=1}^\infty\|T_i(x)\|^2=V\left(\lambda_{\max}\left(\sum_{i=1}^\infty T_i'\circ T_i\right)\right)\cap \E_{H}\subseteq\suppv\left(\sum_{i=1}^\infty T_i' \circ T_i\right).$$ If $v\in \suppv\left(\sum_{i=1}^\infty T_i' \circ T_i\right)$ and 
\begin{equation}\label{loquefalta}
\left\|\left(\sum_{i=1}^\infty T_i' \circ T_i\right)(v)\right\|\leq \sum_{i=1}^\infty\|T_i(v)\|^2,
\end{equation}
then $v\in \argmax_{\|x\|=1} \sum_{i=1}^\infty\|T_i(x)\|^2$.
\end{theorem}

\begin{proof}
First off, since $(T_i)_{i\in\N}$ is a sequence of compact operators, we have that their adjoints are also compact and the composition $T_i’\circ T_i$ is also compact. Since $\mathcal{K}(H)$ is a closed subalgebra of $\Bo(H)$, we have that $\sum_{i=1}^\infty T_i'\circ T_i\in\mathcal{K}(H)$. This means that $$\lambda_{\max}\left(\sum_{i=1}^\infty T_i'\circ T_i\right)=\left\|\sum_{i=1}^\infty T_i'\circ T_i\right\|.$$ In virtue of Lemma \ref{gsv0}, for any $x\in H$ we have that
\begin{eqnarray*}
\sum_{i=1}^\infty\|T_i(x)\|^2 &=& x^*\left(\left(\sum_{i=1}^\infty T_i'\circ T_i\right)(x)\right)\\
&\leq&  \|x^*\|\left\|\sum_{i=1}^\infty T_i'\circ T_i\right\|\|x\|\\
&=& \lambda_{\max}\left(\sum_{i=1}^\infty T_i'\circ T_i\right)\|x\|^2.
\end{eqnarray*}
Therefore $$\max_{\|x\|=1} \sum_{i=1}^\infty\|T_i(x)\|^2 \leq  \lambda_{\max}\left(\sum_{i=1}^\infty T_i'\circ T_i\right).$$ Now, let $w\in V\left(\lambda_{\max}\left(\sum_{i=1}^\infty T_i'\circ T_i\right)\right)\cap \E_{H}$. Then by applying Lemma \ref{gsv0} again
\begin{eqnarray*}
\sum_{i=1}^\infty\|T_i(w)\|^2 &=& w^*\left(\left(\sum_{i=1}^\infty T_i'\circ T_i\right)(w)\right)\\
&=& w^*\left( \lambda_{\max}\left(\sum_{i=1}^\infty T_i'\circ T_i\right)w\right)\\
&=& \lambda_{\max}\left(\sum_{i=1}^\infty T_i'\circ T_i\right).
\end{eqnarray*}
This shows that $$\max_{\|x\|=1} \sum_{i=1}^\infty\|T_i(x)\|^2 = \lambda_{\max}\left(\sum_{i=1}^\infty T_i'\circ T_i\right)$$ and $$V\left(\lambda_{\max}\left(\sum_{i=1}^\infty T_i'\circ T_i\right)\right)\cap \E_{H}\subseteq \argmax_{\|x\|=1} \sum_{i=1}^\infty\|T_i(x)\|^2.$$ Let $\displaystyle{v\in \argmax_{\|x\|=1} \sum_{i=1}^\infty\|T_i(x)\|^2 }$. One the one hand, $$\left\|\frac{\left(\sum_{i=1}^\infty T_i'\circ T_i\right)(v)}{\lambda_{\max}\left(\sum_{i=1}^\infty T_i'\circ T_i\right)}\right\| \leq \frac{\left\|\sum_{i=1}^\infty T_i'\circ T_i\right\|}{\lambda_{\max}\left(\sum_{i=1}^\infty T_i'\circ T_i\right)}= 1.$$ On the other hand and by relying again on Lemma \ref{gsv0}, $$v^*\left(\frac{\left(\sum_{i=1}^\infty T_i'\circ T_i\right)(v)}{\lambda_{\max}\left(\sum_{i=1}^\infty T_i'\circ T_i\right)}\right)=\frac{\sum_{i=1}^\infty\|T_i(v)\|^2}{\lambda_{\max}\left(\sum_{i=1}^\infty T_i'\circ T_i\right)}=1,$$ which implies that $$\left\|\frac{\left(\sum_{i=1}^\infty T_i'\circ T_i\right)(v)}{\lambda_{\max}\left(\sum_{i=1}^\infty T_i'\circ T\right)}\right\| = 1.$$ The strict convexity of $H$ allows us to deduce that $$\frac{\left(\sum_{i=1}^\infty T_i'\circ T_i\right)(v)}{\lambda_{\max}\left(\sum_{i=1}^\infty T_i'\circ T_i\right)}=v,$$ that is, $$\left(\sum_{i=1}^\infty T_i'\circ T_i\right)(v)=\lambda_{\max}\left(\sum_{i=1}^\infty T_i'\circ T_i\right)v$$ and so $v\in V\left(\lambda_{\max}\left(\sum_{i=1}^\infty T_i'\circ T_i\right)\right)\cap \E_{H}$. Next, since $$\lambda_{\max}\left(\sum_{i=1}^\infty T_i'\circ T_i\right)=\left\|\sum_{i=1}^\infty T_i'\circ T_i\right\|,$$ by the observations made at the beginning of this section, we conclude that $$V\left(\lambda_{\max}\left(\sum_{i=1}^\infty T_i'\circ T_i\right)\right)\cap \E_H\subseteq \suppv\left(\sum_{i=1}^\infty T_i'\circ T_i\right).$$ Finally, if $v\in \suppv\left(\sum_{i=1}^\infty T_i' \circ T_i\right)$ and Equation \eqref{loquefalta} holds, then
\begin{eqnarray*}
\max_{\|x\|=1} \sum_{i=1}^\infty\|T_i(x)\|^2 &=& \lambda_{\max}\left(\sum_{i=1}^\infty T_i'\circ T_i\right)\\
&=&\left\|\sum_{i=1}^\infty T_i'\circ T_i\right\|\\
&=&\left\|\left(\sum_{i=1}^\infty T_i' \circ T_i\right)(v)\right\|\\
&\leq& \sum_{i=1}^\infty\|T_i(v)\|^2\\
&\leq& \max_{\|x\|=1} \sum_{i=1}^\infty\|T_i(x)\|^2,
\end{eqnarray*}
which implies that $v\in \argmax_{\|x\|=1} \sum_{i=1}^\infty\|T_i(x)\|^2$.
\end{proof}

In the previous theorem we cannot assure that $$\argmax_{\|x\|=1} \sum_{i=1}^\infty\|T_i(x)\|^2=V\left(\lambda_{\max}\left(\sum_{i=1}^\infty T_i'\circ T_i\right)\right)\cap \E_{H}=\suppv\left(\sum_{i=1}^\infty T_i' \circ T_i\right).$$ The last assertion of Theorem \ref{gsv} states that if $v\in \suppv\left(\sum_{i=1}^\infty T_i' \circ T_i\right)$ and Equation \eqref{loquefalta} holds, then $v\in \argmax_{\|x\|=1} \sum_{i=1}^\infty\|T_i(x)\|^2$. However, in virtue of Lemma \ref{gsv0}, $$\sum_{i=1}^\infty\|T_i(v)\|^2\leq  \left\|\left(\sum_{i=1}^\infty T_i' \circ T_i\right)(v)\right\|,$$ concluding that $$\sum_{i=1}^\infty\|T_i(v)\|^2= \left\|\left(\sum_{i=1}^\infty T_i' \circ T_i\right)(v)\right\|,$$ which in general does not occur.

As a consequence, we cannot assure that the reformulations \eqref{refor0} and \eqref{refor02} are equivalent, but at least we can state that the set of solutions of \eqref{refor0} is contained in the set of solutions of \eqref{refor02}. An schematic summary follows:

$$\eqref{refor0}\Leftrightarrow  \eqref{refor01}\Leftrightarrow  \eqref{refor03}\Rightarrow \eqref{refor02}.$$

\section{Supporting vectors of statistically-normalized operators}

If $x\in\R^m$, then $$\mu_x:=\frac{1}{m}\sum_{i=1}^mx_i$$ and $$\sigma_x^2:=\frac{1}{m}\sum_{i=1}^mx_i^2-\left(\frac{1}{m}\sum_{i=1}^mx_i\right)^2.$$

A vector $x\in\R^m$ is said to be statistically normalized provided that $\mu_x=0$ and $\sigma_x=1$. Note that if $\R^m$ has an statistically normalized vector, then $m\geq 2$.

If $x\in\R^m$ is so that $\sigma_x\neq 0$, then $$x^{st}:=\frac{x-\mu_x{\bf 1}}{\sigma_x}$$ is a statistically normalized vector, where ${\bf 1}$ is the vector whose components are all equal to $1$. Statistically normalized vectors verify special properties.

\begin{remark}\label{normstat}
Let $x,y\in \R^m$ be statistically normalized vectors. Then:
\begin{enumerate}
\item $\|x\|_2^2=m$.
\item $|x\cdot y|\leq m$.
\item $x\cdot y=\frac{\|x+y\|^2_2-2m}{2}$.
\end{enumerate}
\end{remark}

An immediate corollary of Remark \ref{normstat} establishes where to look for statistically normalized vectors of $\R^m$. The place to search within is the sphere of a certain ball centered at $0$ of a subspace of $\R^m$ of codimension $1$. We recall the reader that by $\cdot$ we refer to the Euclidean scalar value and that $\ell_2^m$ stands for $\left(\R^m,\|\bullet\|_2\right)$.

\begin{corollary}\label{snvv}
The set of statistically normalized vectors of $\R^m$ equals $$\E_{\ell_2^m}\left(0,\sqrt{m}\right)\cap \ker({\bf 1}^*).$$
\end{corollary}

\begin{proof}
Simply notice that the set of statistically normalized vectors coincide with the set of solutions of the following nonlinear system:
\begin{equation}\label{snv}
\left\{\begin{array}{l} x_1+\cdots +x_m=0\\ x_1^2+\cdots +x_m^2=m  \end{array}\right.
\end{equation}
\end{proof}

\begin{remark}
The solutions to the system \eqref{snv} can actually be explicitly expressed as follows:
\begin{itemize}
\item If $m=2$, then $\E_{\ell_2^2}\left(0,\sqrt{2}\right)\cap \ker({\bf 1}^*)=\left\{\left(1,-1\right),(-1,1)\right\}$.
\item If $m>2$, then $\E_{\ell_2^m}\left(0,\sqrt{m}\right)\cap \ker({\bf 1}^*)$ is composed of all vectors $x=(x_1,\dots,x_m)$ such that 
$$x_1=\frac{-\sum_{i=3}^mx_i}{2}\mp \sqrt{\frac{m-\sum_{i=3}^mx_i^2}{2}-\frac{\left(\sum_{i=3}^mx_i\right)^2}{4}},$$ $$x_2=\frac{-\sum_{i=3}^mx_i}{2}\pm \sqrt{\frac{m-\sum_{i=3}^mx_i^2}{2}-\frac{\left(\sum_{i=3}^mx_i\right)^2}{4}}$$ and $$2m-2\sum_{i=3}^mx_i^2-\left(\sum_{i=3}^mx_i\right)^2\geq 0.$$
\end{itemize}
\end{remark}

A matrix is said to be statistically normalized if all of its column vectors are statistically normalized. According to Remark \ref{normstat}(1), all the column vectors of an statistically normalized matrix have the same Euclidean norm.

\begin{definition}
Let $X$ be a Banach space with a normalized Schauder basis $(e_n)_{n\in\N}\subseteq \E_X$. A bounded operator $T:X\to X$ is called basic-normalized if $\|T(e_n)\|=1$ for every $n\in\N$.
\end{definition}

Statistically normalized matrices are in fact basic-normalized matrices. Computing the supporting vectors of those matrices of order $m\times 2$ is easier than in the general case.

\subsection{Basic-normalized matrices of $\R^{m\times 2}$}

Even though the proof of the following lemma is nearly trivial, we will include it for the sake of completeness.

\begin{lemma}\label{help}
Let $a,b\in \R$ and consider the function $f(x,y)=a+bxy$ for $(x,y)\in\E_{\ell_2^2}$. Then $\max(f)=a+\frac{|b|}{2}$. In addition:
\begin{enumerate}
\item If $b>0$, then $\argmax(f)=\left\{\left(\frac{\sqrt{2}}{2},\frac{\sqrt{2}}{2}\right),\left(-\frac{\sqrt{2}}{2},-\frac{\sqrt{2}}{2}\right)\right\}$.
\item If $b<0$, then $\argmax(f)=\left\{\left(-\frac{\sqrt{2}}{2},\frac{\sqrt{2}}{2}\right),\left(\frac{\sqrt{2}}{2},-\frac{\sqrt{2}}{2}\right)\right\}$.
\end{enumerate}
\end{lemma}

\begin{proof}
We will consider the case of $b>0$. The other case is similar. We will assume without loss of generality that $y=\sqrt{1-x^2}$ and then we obtain the function of one variable $f(x)=a+bx\sqrt{1-x^2}$ for $x\in[-1,1]$, whose derivative is $$f'(x)=b\sqrt{1-x^2}-\frac{bx^2}{\sqrt{1-x^2}}= \frac{b(1-2x^2)}{\sqrt{1-x^2}}$$ for $x\in(-1,1)$. Note that $f'(x)=0$ means that $x=\pm\frac{\sqrt{2}}{2}$. We spare the rest of the details of the proof to the reader.
\end{proof}

\begin{theorem}\label{teor3.5}
Let $A\in \R^{m\times 2}$ be a matrix whose column vectors $a_1$ and $a_2$ have the same Euclidean norm. Then $\|A\|_2^2=\|a_1\|_2^2+\left|a_1^T\cdot a_2^T\right|$. In addition:
\begin{enumerate}
\item If $a_1^T\cdot a_2^T = 0$, then $\suppv(A)=\E_{\ell_2^2}$.
\item If $a_1^T\cdot a_2^T> 0$, then  $\suppv(A)=\left\{\left(\frac{\sqrt{2}}{2},\frac{\sqrt{2}}{2}\right)^T,\left(-\frac{\sqrt{2}}{2},-\frac{\sqrt{2}}{2}\right)^T\right\}$.
\item If $a_1^T\cdot a_2^T<0$, then $\suppv(A)=\left\{\left(-\frac{\sqrt{2}}{2},\frac{\sqrt{2}}{2}\right)^T,\left(\frac{\sqrt{2}}{2},-\frac{\sqrt{2}}{2}\right)^T\right\}$.
\end{enumerate}
\end{theorem}

\begin{proof}
Let $x^T=(x_1,x_2)\in \E_{\ell_2^2}$ and let us write $A=(a_{ij})$ for $1\leq i\leq m$ and $j=1,2$ where $a_j=(a_{1j}, a_{2j}, \dots a_{mj})^T$ for $j=1,2$. Notice that 
\begin{equation*}
Ax=\left(\begin{array}{c}a_{11}x_1+a_{12}x_2\\ a_{21}x_1+a_{22}x_2\\ \vdots \\ a_{m1}x_1+a_{m2}x_2 \end{array}\right)
\end{equation*}
and
\begin{eqnarray*}
\|Ax\|^2_2&=&\sum_{i=1}^m  (a_{i1}x_1+a_{i2}x_2)^2\\
&=&x_1^2\sum_{i=1}^m  a_{i1}^2 + x_2^2\sum_{i=1}^ma_{i2}^2+ 2x_1x_2\sum_{i=1}^ma_{i1}a_{i2}\\
&=& (x_1^2+x_2^2)\sum_{i=1}^m  a_{i1}^2+2x_1x_2\sum_{i=1}^ma_{i1}a_{i2}\\
&=& \|a_1\|_2^2+ 2x_1x_2a_1^T\cdot a_2^T.
\end{eqnarray*}
At this stage, it only suffices to apply Lemma \ref{help}.
\end{proof}

\subsection{Basic-normalized matrices of $\R^{m\times n}$}

The general case of basic-normalized matrices is more complex than the previous one as we will show next. Except for particular cases, it is not worth solving the general case in a straight forward way as before.

\begin{lemma}\label{help2}
Let $A\in \R^{m\times n}$ be a matrix whose column vectors $a_1,\dots, a_n$ have the same Euclidean norm. Then $$\|Ax\|_2^2=\|a_1\|_2^2 + 2\sum_{j\neq k}x_jx_ka_j^T\cdot a_k^T$$ for every $x^T=(x_1,\dots,x_n)\in\E_{\ell_2^n}$.
\end{lemma}

\begin{proof}
Let $x^T=(x_1,x_2,\dots,x_n)\in \E_{\ell_2^n}$ and let us write $A=(a_{ij})$ for $1\leq i\leq m$ and $1\leq j\leq n$ where $a_j=(a_{11}, a_{21}, \dots a_{mj})^T$ for $j=1,2,\dots,n$. Notice that 
\begin{equation*}
Ax=\left(\begin{array}{c}a_{11}x_1+a_{12}x_2+\cdots +a_{1n}x_n\\ a_{21}x_1+a_{22}x_2+\cdots+a_{2n}x_n\\ \vdots \\ a_{m1}x_1+a_{m2}x_2+\cdots+a_{mn}x_n \end{array}\right)
\end{equation*}
and
\begin{eqnarray*}
\|Ax\|^2_2&=& \sum_{i=1}^m\left(\sum_{j=1}^na_{ij}x_j\right)^2\\
&=&\sum_{i=1}^m\sum_{j=1}^na_{ij}^2x_j^2 + 2\sum_{i=1}^m\sum_{j\neq k}a_{ij}x_ja_{ik}x_k\\
&=&\sum_{j=1}^nx_j^2\sum_{i=1}^ma_{ij}^2 + 2\sum_{j\neq k}x_jx_k\sum_{i=1}^ma_{ij}a_{ik}\\
&=&\|a_1\|_2^2 + 2\sum_{j\neq k}x_jx_ka_j^T\cdot a_k^T.
\end{eqnarray*}
\end{proof}

As we can see from the proof of the previous theorem, in order to find the supporting vectors of a matrix $A\in \R^{m\times n}$ whose column vectors $a_1,\dots, a_n$ have the same Euclidean norm, all we need to do is maximize a function of the form $f(x_1,\dots,x_n)=a+\sum_{j\neq k}c_{jk}x_jx_k$ for $(x_1,\dots,x_k)\in\E_{\ell_2^n}$ where $a>0$ and $c_{jk}\in \R$ for $j\neq k\in\{1,\dots,n\}$ satisfying that $c_{jk}=c_{kj}$. If we apply the Lagrange multiplier technique to $f$, then we obtain the function $$\mathcal{L}(x_1,\dots,x_n,\lambda):=a+\sum_{j\neq k}c_{jk}x_jx_k- \lambda\left(1-\sum_{j=1}^n x_j^2\right)$$ whose partial derivatives are $$\frac{\partial \mathcal{L}}{\partial x_j} = \sum_{k\neq j} c_{jk}x_k + 2\lambda x_j$$ for all $j\in\{1,\dots, n\}$. The critical points of $\mathcal{L}$ verify the following nonlinear system \begin{equation}\label{crit}\left\{\begin{array}{l} \left(\begin{array}{cccc} 2\lambda & c_{12} & \cdots & c_{1n}\\ c_{21} & 2\lambda & \cdots & c_{2n} \\ \vdots & \vdots & \ddots & \vdots \\ c_{n1} & c_{n2} & \cdots & 2\lambda \end{array}\right)\left(\begin{array}{c}x_1\\ x_2\\ \vdots \\ x_n\end{array}\right)=\left(\begin{array}{c}0\\ 0\\ \vdots \\ 0\end{array}\right)\\ \\ x_1^2+x_2^2+\cdots + x_n^2=1 \end{array}\right.\end{equation} Notice that if the coefficient matrix above is invertible, then the previous nonlinear system has no solution, so $\lambda$ must be found to verify that $$\left|\begin{array}{cccc}2\lambda & c_{12} & \cdots & c_{1n}\\ c_{21} & 2\lambda & \cdots & c_{2n} \\ \vdots & \vdots & \ddots & \vdots \\ c_{n1} & c_{n2} & \cdots & 2\lambda \end{array}\right|=0.$$ Bear in mind that the above matrix is symmetric. If the coefficients $c_{jk}$ verify special properties, then we can throw some light into the critical point of $\mathcal{L}$. The following lemma deals with the case where all the $c_{jk}$’s are equal to each other.

\begin{lemma}\label{partcase}
Let $f(x_1,\dots,x_n)=a+c\sum_{j\neq k}x_jx_k$ for $(x_1,\dots,x_k)\in\E_{\ell_2^n}$ where $a>0$ and $c\in \R\setminus\{0\}$. Then the norm-$1$ multiples of the statistically normalized vectors of $\R^n$ are critical points of the Lagrangian function $$\mathcal{L}(x_1,\dots,x_n,\lambda):=a+c\sum_{j\neq k}x_jx_k- \lambda\left(1-\sum_{j=1}^n x_j^2\right).$$
\end{lemma}

\begin{proof}
We just need to go back to the nonlinear system \eqref{crit} and adapt it to our hypothesis to obtain
\begin{equation}\label{crit2}
\left\{\begin{array}{l} \left(\begin{array}{cccc} 2\lambda & c & \cdots & c\\ c & 2\lambda & \cdots & c \\ \vdots & \vdots & \ddots & \vdots \\ c & c & \cdots & 2\lambda \end{array}\right)\left(\begin{array}{c}x_1\\ x_2\\ \vdots \\ x_n\end{array}\right)=\left(\begin{array}{c}0\\ 0\\ \vdots \\ 0\end{array}\right)\\ \\ x_1^2+x_2^2+\cdots + x_n^2=1 \end{array}\right.
\end{equation}
and to realize that $\lambda= \frac{c}{2}$ turns the previous system into $$
\left\{\begin{array}{l} x_1+\cdots +x_n=0\\ x_1^2+\cdots +x_n^2=1  \end{array}\right.$$ which gives us the desired result (see Corollary \ref{snvv} together with Equation \eqref{snv}).
\end{proof}

The following remark shows that, in the settings of Lemma \ref{partcase}, not all the critical points of the Lagrangian function are norm-$1$ multiples of the statistically normalized vectors of $\R^n$. 

\begin{remark}
Observe that $$\left|\begin{array}{ccc}2\lambda & c & c\\ c & 2\lambda & c \\ c & c & 2\lambda \end{array}\right|=2(2\lambda-c)^2(\lambda+c)$$ so $\lambda=-c$ also provides critical points whose mean is not null. Indeed, if $\lambda =-c$, then system \eqref{crit2} turns into $$\left\{\begin{array}{l}\left(\begin{array}{ccc} -2c&c&c\\ c&-2c&c\\ c&c&-2c \end{array}\right)\left(\begin{array}{c} x_1\\x_2\\x_3 \end{array}\right)=\left(\begin{array}{c}0\\ 0\\ 0 \end{array}\right) \\ \\ x_1^2+x_2^2+x_3^2=1\end{array}\right.$$ whose set of solutions is $$\left\{\left(\frac{ 1}{\sqrt{3}},\frac{ 1}{\sqrt{3}},\frac{ 1}{\sqrt{3}}\right),\left(\frac{-1}{\sqrt{3}},\frac{-1}{\sqrt{3}},\frac{- 1}{\sqrt{3}}\right)\right\}.$$
\end{remark}

After this remark we can state without proof the following lemma.

\begin{lemma}\label{partcase2}
Let $f(x_1,x_2,x_3)=a+c\sum_{j\neq k}x_jx_k$ for $(x_1,x_2,x_3)\in\E_{\ell_2^3}$ where $a>0$ and $c\in \R\setminus\{0\}$. Then the critical points of the Lagrangian function $$\mathcal{L}(x_1,x_2,x_3,\lambda):=a+c\sum_{j\neq k}x_jx_k- \lambda\left(1-\sum_{j=1}^3 x_j^2\right)$$ consist exactly of the norm-$1$ multiples of the statistically normalized vectors of $\R^3$ together with $\left\{\left(\frac{ 1}{\sqrt{3}},\frac{ 1}{\sqrt{3}},\frac{ 1}{\sqrt{3}}\right),\left(\frac{-1}{\sqrt{3}},\frac{-1}{\sqrt{3}},\frac{- 1}{\sqrt{3}}\right)\right\}.$
\end{lemma}

To summarize, except for particular cases where the coefficients $c_{jk}$ verify special properties like in Lemma \ref{partcase}, in order to find the supporting vectors of a matrix $A$ like in Lemma \ref{help2} it is more efficient to rely on Theorem \ref{gsv}.

\section{Applications of the proposed model (Theorem \ref{gsv})}

\subsection{Pure state that jointly maximizes the modulus of observable magnitudes in a quantum mechanical system, with special attention to the probability density operator}\label{quantum}

Let us begin by recalling the first two postulates of Quantum Mechanics (see \cite{SN}). We will not follow the classical Quantum Mechanics notation (the Bra-Ket notation), but the classical Functional Analysis notation to keep consistency with the notation in the rest of the paper:

\begin{enumerate}

\item The First Postulate of Quantum Mechanics establishes that to every quantum mechanical system an infinite dimensional separable complex Hilbert space $H$ corresponds. A pure state of this system in a fixed instant of time $t$ is represented by a unit ray $\E_{\C}x$ with $\|x\|=1$. An element of the previous ray is called a state vector or a ket.

\item The Second Postulate of Quantum Mechanics establishes that every observable magnitude of the quantum mechanical system $H$ is represented by a selfadjoint linear operator $T: H \to H$. This correspondence between observable magnitudes and selfadjoint linear operators is not in general bijective, that is, not all selfadjoint linear operators represent an observable magnitude. The existence of observable magnitudes represented by a selfadjoint {\bf unbounded} operators implies that the Hilbert space representing the quantum mechanical system is infinite dimensional, since every linear operator on a finite dimensional Banach space is compact and thus bounded. If an observable magnitude is represented by a selfadjoint bounded operator $T: H\to H$, then $\|T\|$ measures the intensity of the observable magnitude. Since $T$ is sefladjoint, the residual spectrum of $T$, $\sigma_r(T)$, is empty and thus the spectrum of $T$ is the disjoint union of the point spectrum and the continuous spectrum, $\sigma(T)=\sigma_p(T)\cup\sigma_c(T)$. Also, $\sigma(T)\subseteq \R$ and if $\lambda \neq \gamma\in \sigma_p(T)$, then $V(\lambda)\subseteq V(\gamma)^{\perp}$. By $\left(e_{(\lambda,n)}\right)_{n\in I(\lambda)}$ we denote an orthonormal basis of $V(\lambda)$ and $I(\lambda):=\dim\left(V(\lambda)\right)$. Thus $\left(e_{(\lambda,n)}\right)_{(\lambda,n)\in\sigma_p(T)\times I(\lambda)}$ is an orthonormal system in $H$.




\end{enumerate}

A probability density matrix represents a partial state of knowledge of a (finite-dimensional) system (see \cite[Section 6]{Grandy2004}): $$\rho(\bullet)=\sum_{i=1}^n w_i(\bullet|\psi_i) \psi_i.$$ Based on that information we conclude that with probability $w_i$ the system
may be in a pure state $\psi_i$.

For quantum systems represented by infinite dimensional complex separable Hilbert spaces (for instance, those with unbounded observable magnitudes), the probability density matrix is in fact an operator, which we describe next.

Let $H$ be an infinite dimensional separable complex Hilbert space representing a quantum mechanical system. In the first place, if the system is in a mixed state given by the following states $(x_n)_{n\in\N}\subseteq \E_H$, where $\E_{\C}x_n\cap \E_{\C}x_m=\varnothing$ if $n\neq m$, then the probability density operator is given by 
\begin{equation}\label{prodensoper}
\begin{array}{rrcl}
    D: & H &\to & H \\
     & x &\mapsto & \displaystyle{ D(x):= \sum_{n=1}^\infty \rho_n x_n^*(x)x_n=\sum_{n=1}^\infty \rho_n (x|x_n)x_n}
\end{array}
\end{equation}
where $\sum_{n=1}^\infty \rho_n$ is a convex series, that is, $\rho_n\geq 0$ for all  $n\in\N$ and $\sum_{n=1}^\infty \rho_n=1$. Note that $\rho_n$ indicates the probability that the system $H$ is at the state $x_n$. The probability density operator $D$ clearly satisfies the following properties:
\begin{itemize}
    \item $D$ is clearly bounded. Indeed, for every $x\in H$, $$\|D(x)\|\leq \sum_{n=1}^\infty \rho_n |(x|x_n)|\|x_n\|\leq \sum_{n=1}^\infty \rho_n\|x\|\|x_n\|^2=\|x\|\sum_{n=1}^\infty \rho_n=\|x\|,$$ which implies that $D$ is continuous and $\|D\|\leq 1$.
    \item $D$ is clearly selfadjoint. Indeed, for every $x,y\in H$
    \begin{eqnarray*}
    (D(x)|y)&=&\left(\sum_{n=1}^\infty \rho_n (x|x_n)x_n\Big |y\right)\\
    &=&\sum_{n=1}^\infty \rho_n (x|x_n)(x_n|y)\\
    &=&\sum_{n=1}^\infty \rho_n \overline{(y|x_n)}(x|x_n)\\
    &=&\left(x\Big|\sum_{n=1}^\infty \rho_n (y|x_n)x_n\right)\\
    &=& (x|D(y)).
    \end{eqnarray*}
    \item $D$ is clearly compact. Indeed, $D$ is compact because it can approximated in the operator norm of $\Bo(H)$ by the sequence $(D_k)_{k\in\N}$ of finite-rank selfadjoint operators, where
    \begin{equation*}
\begin{array}{rrcl}
    D_k: & H &\to & \spa\{x_1,\dots,x_k\} \\
     & x &\mapsto & \displaystyle{ D(x):= \sum_{n=1}^k \rho_n x_n^*(x)x_n=\sum_{n=1}^k \rho_n (x|x_n)x_n}.
\end{array}
\end{equation*} Indeed, for every $x\in \B_H$ we have that $$\|D(x)-D_k(x)\|=\left\|\sum_{n=k+1}^\infty \rho_n x_n^*(x)x_n\right\|\leq \sum_{n=k+1}^\infty \rho_n\to 0$$ as $k\to \infty$ because $\sum_{n=k+1}^\infty \rho_n$ is the rest of a convergent series.

\item $D$ is clearly positive, in fact, $D\geq D^2\geq 0$. Indeed, $$(D(x)|x)=\sum_{n=1}^\infty \rho_n (x|x_n)(x_n|x)=\sum_{n=1}^\infty\rho_n |(x|x_n)|^2\geq 0 $$ for all $x\in X$. Note that $D^2$ is a normal operator (positive and selfadjoint). Let us prove now that $D\geq D^2$. Fix an arbitrary $x\in H$. On the one hand, $$ \left(D^2(x)|x\right) = \sum_{n=1}^\infty \rho_n(D(x)|x_n)(x_n|x)= \sum_{n=1}^\infty \rho_n\sum_{k=1}^\infty \rho_k (x|x_k)(x_k|x_n)(x_n|x).$$ On the other hand, since $D$ is selfadjoint, $\left(D^2(x)|x\right)=(D(x)|D(x))=\|D(x)\|^2$, so $\left(D^2(x)|x\right)$ is real and positive. Then, H\"older’s inequality allows that
\begin{eqnarray*}
\left(D^2(x)|x\right)&=& \Re\left(\left(D^2(x)|x\right)\right)\\
&=& \sum_{n=1}^\infty \rho_n\sum_{k=1}^\infty \rho_k\Re\left((x|x_k)(x_k|x_n)(x_n|x)\right)\\
&\leq & \sum_{n=1}^\infty\rho_n\sum_{k=1}^\infty \rho_k\left|(x|x_k)(x_k|x_n)(x_n|x)\right|\\
&\leq  & \sum_{n=1}^\infty\rho_n|(x_n|x)|\sum_{k=1}^\infty \rho_k|(x|x_k)|\\
&=&\left(\sum_{n=1}^\infty \rho_n|(x_n|x)|\right)^2\\
&=&\left(\sum_{n=1}^\infty \sqrt{\rho_n}\sqrt{\rho_n}|(x_n|x)|\right)^2\\
&\leq  & \left(\sum_{n=1}^\infty \rho_n\right)\left(\sum_{n=1}^\infty \rho_n|(x_n|x)|^2\right)\\
&=& (D(x)|x).
\end{eqnarray*}

\item $\mathrm{tr}(D)=1$. Indeed, let $(e_k)_{k\in\N}$ be an orthonormal basis of $H$, then
\begin{eqnarray*}
    \mathrm{tr}(D)&=&\sum_{k=1}^\infty e_k^*(D(e_k))=\sum_{k=1}^\infty (D(e_k)|e_k)\\
    &=&\sum_{k=1}^\infty \left(\sum_{n=1}^\infty \rho_n (e_k|x_n)x_n\Big |e_k\right)=\sum_{k=1}^\infty \sum_{n=1}^\infty \rho_n (e_k|x_n)\left(x_n |e_k\right)\\
    &=&\sum_{k=1}^\infty \sum_{n=1}^\infty \rho_n \overline{(x_n |e_k)}(x_n |e_k)=\sum_{k=1}^\infty \sum_{n=1}^\infty \rho_n \left|(x_n |e_k)\right|^2\\
     &=&\sum_{n=1}^\infty \rho_n\sum_{k=1}^\infty  \left|(x_n |e_k)\right|^2=\sum_{n=1}^\infty \rho_n\|x_n \|^2=\sum_{n=1}^\infty \rho_n= 1.
\end{eqnarray*}

\item If $(x_n)_{n\in\N}$ is an orthonormal basis, then $\|D\|=\displaystyle{\left\|\left(\rho_n\right)_{n\in\N}\right\|_\infty}$. Indeed, let $x\in\E_H$ and write $x=\sum_{n=1}^\infty (x|x_n)x_n$, where $1=\sum_{n=1}^\infty |(x|x_n)|^2$. Observe that 
\begin{eqnarray*}
\|D(x)\|&=&\sqrt{\sum_{n=1}^\infty \rho_n^2|(x|x_n)|^2}\\
&\leq& \sqrt{\sum_{n=1}^\infty \left\|\left(\rho_n^2\right)_{n\in\N}\right\|_\infty|(x|x_n)|^2}\\
&=&\left\|\left(\rho_n\right)_{n\in\N}\right\|_\infty\sqrt{\sum_{n=1}^\infty |(x|x_n)|^2}\\
&=& \left\|\left(\rho_n\right)_{n\in\N}\right\|_\infty.
\end{eqnarray*}
This shows that $\|D\|\leq \displaystyle{\left\|\left(\rho_n\right)_{n\in\N}\right\|_\infty}$. In order to see that $\|D\|=\displaystyle{\left\|\left(\rho_n\right)_{n\in\N}\right\|_\infty}$, we just need to choose $n_0\in\N$ such that $\displaystyle{\rho_{n_0}=\left\|\left(\rho_n\right)_{n\in\N}\right\|_\infty=\max_{n\in\N}\rho_n}$ and realize that $D(x_{n_0})=\rho_{n_0}x_{n_0}$, therefore,
\begin{equation}\label{suppvec}
\|D(x_0)\|=\rho_{n_0}=\left\|\left(\rho_n\right)_{n\in\N}\right\|_\infty.
\end{equation}
This shows that $\|D\|=\displaystyle{\left\|\left(\rho_n\right)_{n\in\N}\right\|_\infty}$. According to Theorem \ref{gsv}, since $D$ is selfadjoint and compact, then $$\argmax_{\|x\|=1}\|D(x)\|^2=V\left(\lambda_{\max}(D^2)\right)\cap\E_H.$$ Notice that $$\argmax_{\|x\|=1}\|D(x)\|^2= \argmax_{\|x\|=1}\|D(x)\|=\suppv(D).$$ A supporting vector of $D$ is $x_{n_0}$ in virtue of Equation \eqref{suppvec}, which is precisely the state with the highest probability of the system.
\end{itemize}

Here we just found a first application of our main result (Theorem \ref{gsv}). A concrete example of a probability density operator can be given in the infinite dimensional separable complex Hilbert space $$H:=\ell_2=\ell_2(\N,\C):=\left\{(\alpha_k)_{k\in\N}\in \C^{\N}:\left\|(\alpha_k)_{k\in\N}\right\|_2:=\left(\sum_{k=1}^\infty |\alpha_k|^2\right)^{\frac{1}{2}}<\infty\right\}$$ with states given by the canonical basis $(e_n)_{n\in\N}$ and probabilities $\rho_n=\frac{1}{2^n}$ for every $n\in\N$. In this situation, the probability density operator looks like 
\begin{equation}\label{prodensoper2}
\begin{array}{rrcl}
    D: & \ell_2 &\to & \ell_2 \\
     & (\alpha_k)_{k\in\N} &\mapsto & \displaystyle{ D\left((\alpha_k)_{k\in\N}\right):= \left(\frac{\alpha_k}{2^k}\right)_{k\in\N}}.
\end{array}
\end{equation}
According to above, $\|D\|=\frac{1}{2}$ and a supporting vector of $D$ is $e_1$, which is the state with highest probability.

More generally, let $(T_n)_{n\in\N}\subset \Bo(H)$ be a sequence of selfadjoint compact operators representing observable magnitudes. For every $n\in \N$, $\|T_n\|$ represents the modulus of the magnitude $T_n$. The pure state of $H$ that jointly maximizes the moduli of the previous magnitudes is $$\argmax_{\|x\|=1}\sum_{n=1}^\infty \|T_n(x)\|^2.$$ According to our Theorem \ref{gsv}, such pure state can be found in $$V\left(\lambda_{\max}\left(\sum_{n=1}^\infty T_n^2\right)\right)\cap \E_{H}.$$

\subsection{Optimal TMS coils}

Transcranial Magnetic Stimulation (TMS) is a noninvasive technique to stimulate the brain, which is applied to psychiatric and medical conditions, such as major depressive disorder, schizophrenia, bipolar depression, post-traumatic, stress disorder and obsessive-compulsive disorder, amongst others \cite{key:wassermann}.
The development of TMS is being restricted by technical limitations, such as the undesired stimulation in non-target cortex regions.  

Over the years, there have been new TMS stimulator design methods such as \cite{KOP,mente1,KOP2,CSGPGRH,Gomez,Wang}. In all these approaches, coil design problem is eventually posed as a convex optimization, where the constant search for new coil features and improved performance has highlighted the need of employing more versatile optimization techniques capable of dealing with the new requirements.

In TMS, strong current pulses driven through a coil are used to induce an electric field stimulating neurons in the cortex.

The goal in TMS coil design is to find optimal positions for the multiple windings of coils (or equivalently the  current density) so as to produce fields with the desired spatial characteristics and properties \cite{mente1,CSGPGRH} (high focality, field penetration depth, low inductance, low heat dissipation, etc.). This design problem has been frequently posed as a constraint optimization problem. The idea is to solve the resulting optimization problem in an exact manner justified by abstract mathematical proofs unlike previous TMS coils designed by means of heuristic methods not proven to be convergent to the optimal solution \cite{mente1,key:wassermann,key:romei}.

In this work, in order to illustrate an application of the maths above, we are going to tackle the design of a novel TMS coil capable of producing a maximal stimulation in the occipital lobe while dissipating a minimum power. The modulation of the activity in this region in the cerebral cortex has an enormous interest as it is involved in many brain functions as those related with process of visual information \cite{key:romei}. 

\begin{figure}[H]
\begin{center}
\subfigure[] {\label{fig:setup2}\includegraphics[width=0.8\textwidth]{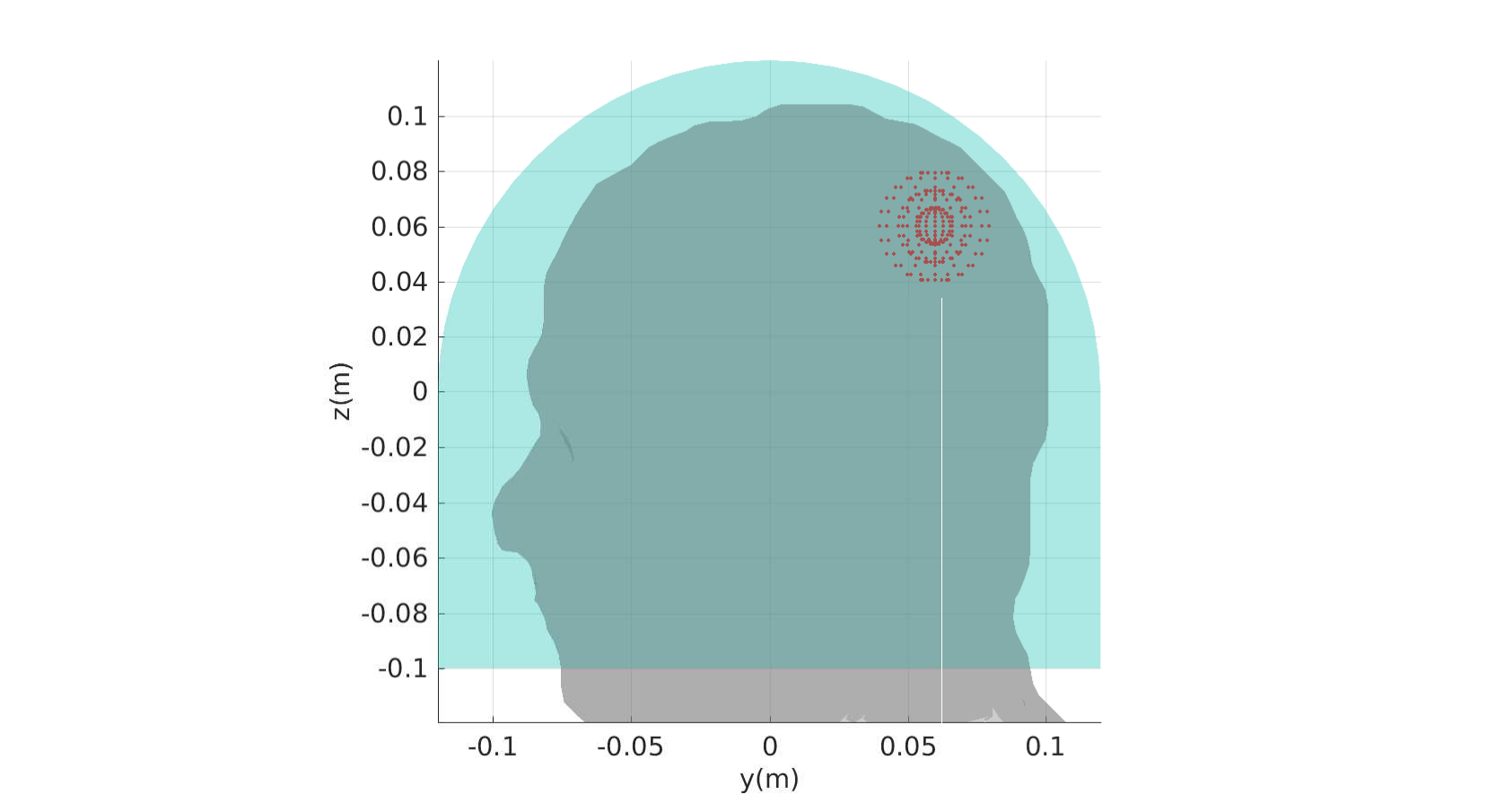}} %
\subfigure[] {\label{fig:setup3}\includegraphics[width=1.0\textwidth]{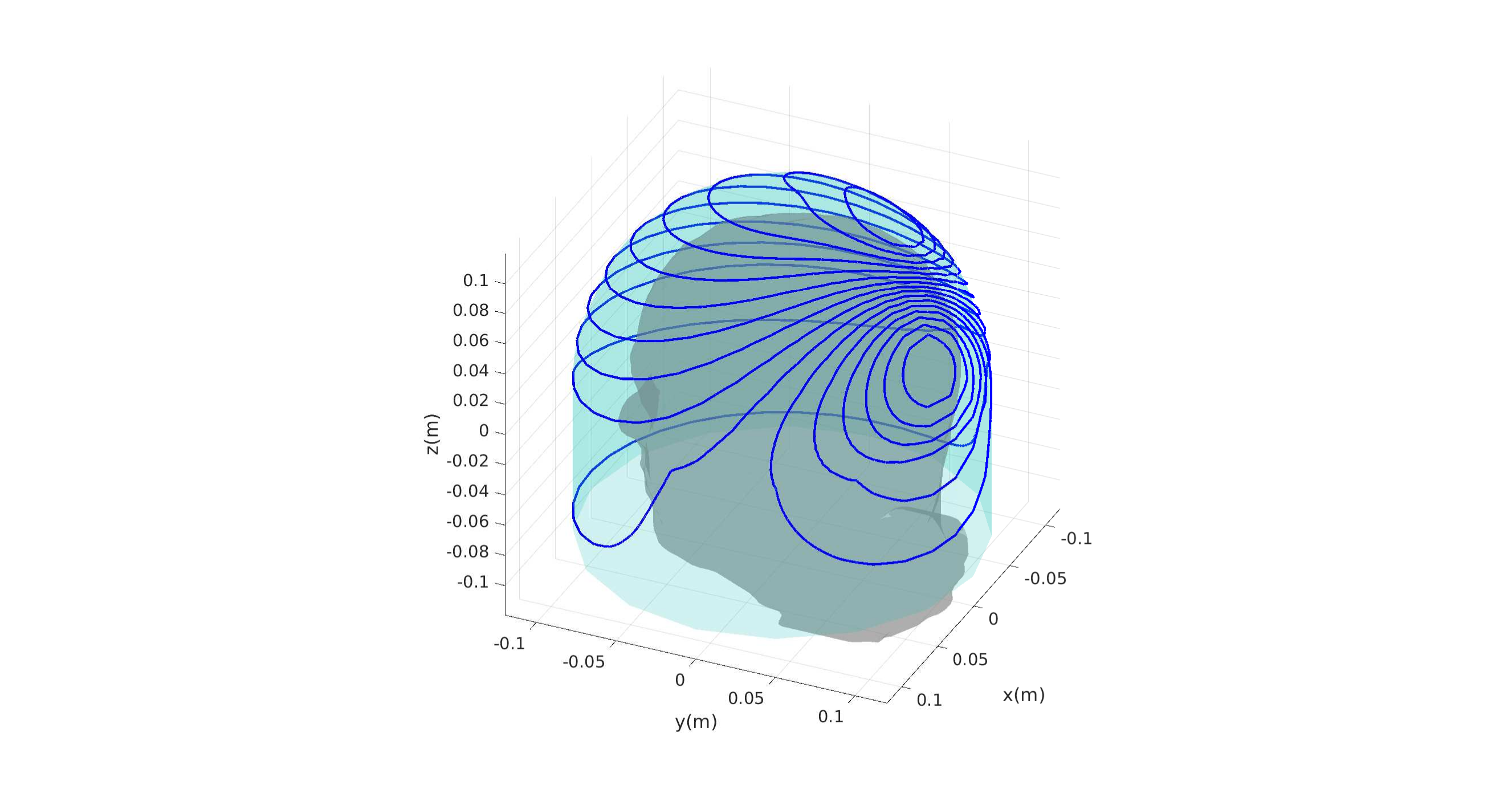}} %
\end{center}
\caption{\emph{(a) Schematic diagram showing the hemispherical conducting surface with a cylindrical extension along with the region of interest where the stimulation is desired to be maximal, which has been included in a illustrative human head model for sake of clarity. (b) Wirepaths with 16 turns of the TMS coil solution of design problem in Equation \eqref{ffm}.}}
\end{figure}

To this end a new coil geometry consisting of a hemispherical surface of radius $0.12m$ with a cylindrical extension of height $0.1m$ has been considered, in which we want to find the electric current that induces maximum $E_x$, $E_y$ and $E_z$ fields in a target region in the occipital lobe formed from a $2cm$ radius spherical distribution of $H=400$ points, of which center is at coordinates $(0cm, 6cm, 6cm)$ as shown in figure \ref{fig:setup2}.

This particular problem can be formulated as

\begin{equation} \label{ffm}
\left\{
\begin{array}{l}
\max  \left\| E_z \psi \right\| _2 \\
\max  \left\| E_y \psi \right\| _2 \\
\max  \left\| E_x \psi \right\| _2 \\
\min   \psi^T R \psi
\end{array}
\right.
\end{equation}
where $N$ is the number of nodes employed to mesh the conducting surface \cite{CSGPGRH}, with $N>H$, $E_x,E_y, E_z \in \mathbb{R}^{H \times N}$, $R \in \mathbb{R}^{N \times N}$ is the resistance and $\psi$ is the desired stream function.

We proceed now to reformulate the multiobjective optimization problem given in Equation \eqref{ffm}. First, we apply the Cholesky decomposition to $R$ to obtain $R=C^TC$ so we have that $\psi^TR\psi =(C\psi)^T(C\psi)=\|C\psi\|_2^2$ and we obtain
 \begin{equation} \label{ffm1}
\left\{
\begin{array}{l}
\max  \left\| E_z \psi \right\| _2 \\
\max  \left\| E_y \psi \right\| _2 \\
\max  \left\| E_x \psi \right\| _2 \\
\min  \|C\psi\|_2^2  
\end{array}
\right.
\end{equation}

Next, by taking into consideration that the square root is a strictly increasing function on $[0,\infty)$, we can apply Equation \eqref{cv2} to obtain 
\begin{equation} \label{ffm2}
\left\{
\begin{array}{l}
\max  \left\| E_z \psi \right\| _2 \\
\max  \left\| E_y \psi \right\| _2 \\
\max  \left\| E_x \psi \right\| _2 \\
\min  \|C\psi\|_2
\end{array}
\right.
\end{equation}

Now, since $C$ is an invertible square matrix, by setting $\varphi=C\psi$ we obtain

\begin{equation} \label{ffm3}
\left\{
\begin{array}{l}
\max  \left\| \left(E_zC^{-1}\right) \varphi \right\| _2 \\
\max  \left\| \left(E_yC^{-1}\right) \varphi \right\| _2 \\
\max  \left\| \left(E_xC^{-1}\right) \varphi \right\| _2 \\
\min  \|\varphi\|_2
\end{array}
\right.
\end{equation}

Problem \eqref{ffm3} is of the form \eqref{original}. In accordance with Section \ref{orgsec}, Problem \eqref{ffm3} is reformulated as follows in order not to fall out of the Hilbert space setting:
\begin{equation} \label{ffm6}
\left\{
\begin{array}{l}
\max  \left\| (E_zC^{-1})\varphi \right\|^2 _2 + \left\| (E_yC^{-1})\varphi \right\|^2 _2 + \left\| (E_xC^{-1})\varphi \right\|^2 _2 \\
\|\varphi\|_2  = 1
\end{array}
\right.
\end{equation}
Problem \eqref{ffm6} can be solved with Theorem \ref{gsv}. Once we find a solution $\varphi$ of Problem \eqref{ffm6}, we have that $\psi=C^{-1}\varphi$ is the desired stream function.

\begin{figure}[H]
\begin{center}
\subfigure[] {\label{fig:clemente.JPG}\includegraphics[width=0.85\textwidth]{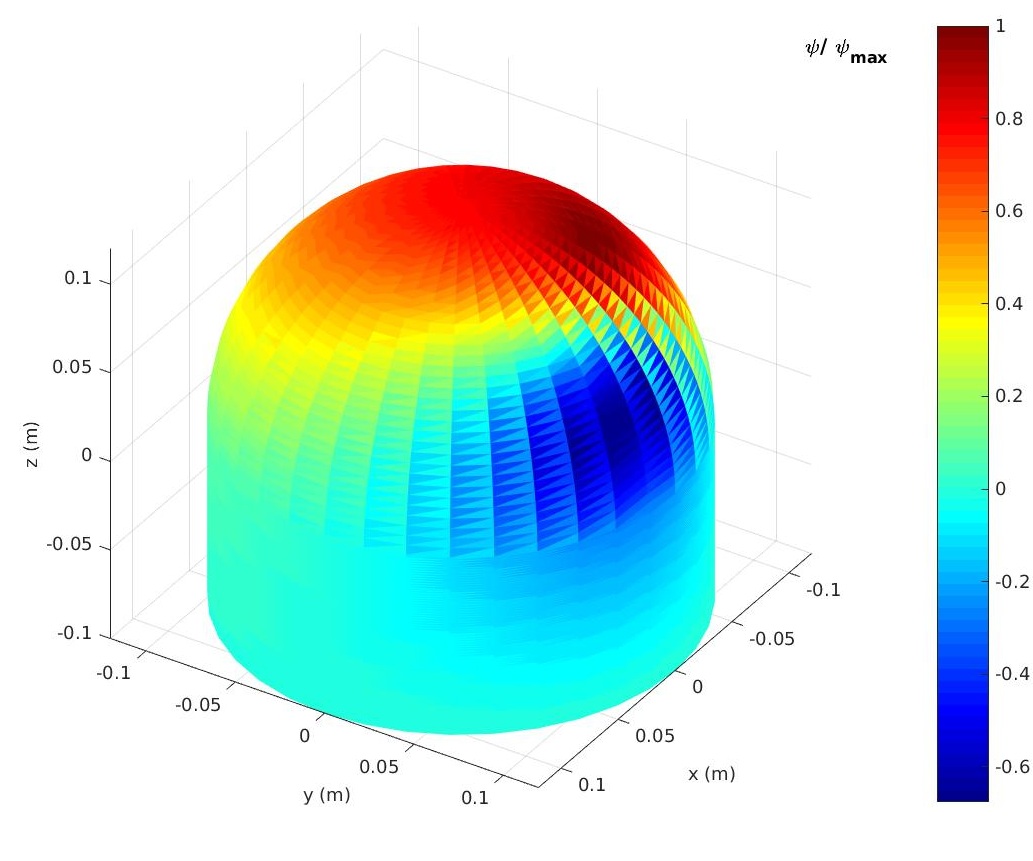}} %
\end{center}
\caption{\emph{(a) Colormap of the normalised optimal stream function $\dfrac{\psi}{ \psi_{\max}}$  over the coil surface.}}
\end{figure}

    The wire-paths of the solution to the design problem in Equation \eqref{ffm} is shown in Figure \ref{fig:setup3}, it is a two lobed TMS coil, where it can be seen that the winding density is more concentrated over the region of stimulation. Finally, Figure \ref{fig:clemente.JPG} shows the colormap of the normalized optimal stream function over the coil surface.

\subsection{Optimal location using statistics multivariate analysis}\label{secc_alb}

In multivariate statistical analysis, the use of different techniques such as Principal Component Analysis (PCA) for modelling the variance between variables, has been largely considered \cite{alber1} with up-to-date applications such as big data and digital image analysis \cite{alber2, alber3}. Here we present a practical application of supporting vectors focus on the optimal location of a private educational academy between the different provinces of Spain.

We are going to consider 3 different variables for this study: proportion of population with some kind of education (basic or high) ($m_1$), with higher education ($m_2$) and with higher education between 25 and 34 years ($m_3$). These data have been obtained from the 2009 social indicators database of the National Statistic Office (in Spain INE. \url{http://www.ine.es/daco/daco42/sociales09/sociales.htm}).

To find the optimal location we consider the province where the values of $m_1$, $m_2$ and $m_3$ are maximum. Due to the different units of the variables, we standardize them ($\mu=0, \sigma=1$) obtaining $m_1^{st}$, $m_2^{st}$ and $m_3^{st}$. Here we have a typical multiobjective problem that can be formulated as follow:

\begin{equation} \label{alb1}
\left\{
\begin{array}{l}
\max  \left\| M x\right\|^2_2 \\
\|x\|_2^2=1
\end{array}
\right.
\end{equation}
where $M$ is the matrix composed by the three variables standardized considered in our study $M=[m_1^{st}|m_2^{st}|m_3^{st}]$ and $x$ is the supporting vector (in this case, the first component in PCA as shown in Figure 3), with $N=52$ locations, $H=3$ standardized variables, $x \in \mathbb{R}^{H \times 1}$ and $m_1^{st}, m_2^{st}, m_3^{st} \in \mathbb{R}^{N \times 1}$. We observe that this optimization problem is similar to Equation \eqref{ffm6} and thus Theorem \ref{gsv} applies to it.

Moreover, the solution of this multiobjective problem let us to sort the sites considering the three variables of this study (Figure 4A). The higher value of $Mx$ indicate the best place to locate the private educational academy because it has the higher proportions of population which can be considered as target. Likewise, we can represent them in a geographic projection (Figure 4B). Although this is an example, it can be extended to more variables achieving a better model and results.

\begin{figure}[H]
\begin{center}
\label{fig:alb1png}\includegraphics[width=1.0\textwidth]{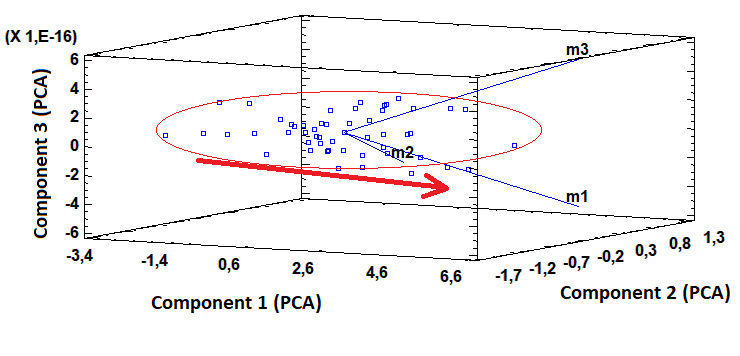} %
\end{center}
\caption{\emph{Reference frame defined by the three principal components with the variables of Section \ref{secc_alb}. In this case, with the variables directly correlated, a supporting vector $x$ indicates the direction of the first component (maximum of the three variables, red arrow). In blue squares, the locations considered in this paper. In blue lines, the standardized variables $m_1^{st}, m_2^{st}, m_3^{st}$ represented in PCA reference frame.}}
\end{figure}

\begin{figure}[H]
\begin{center}
\label{fig:Mxpng}\includegraphics[width=1\textwidth]{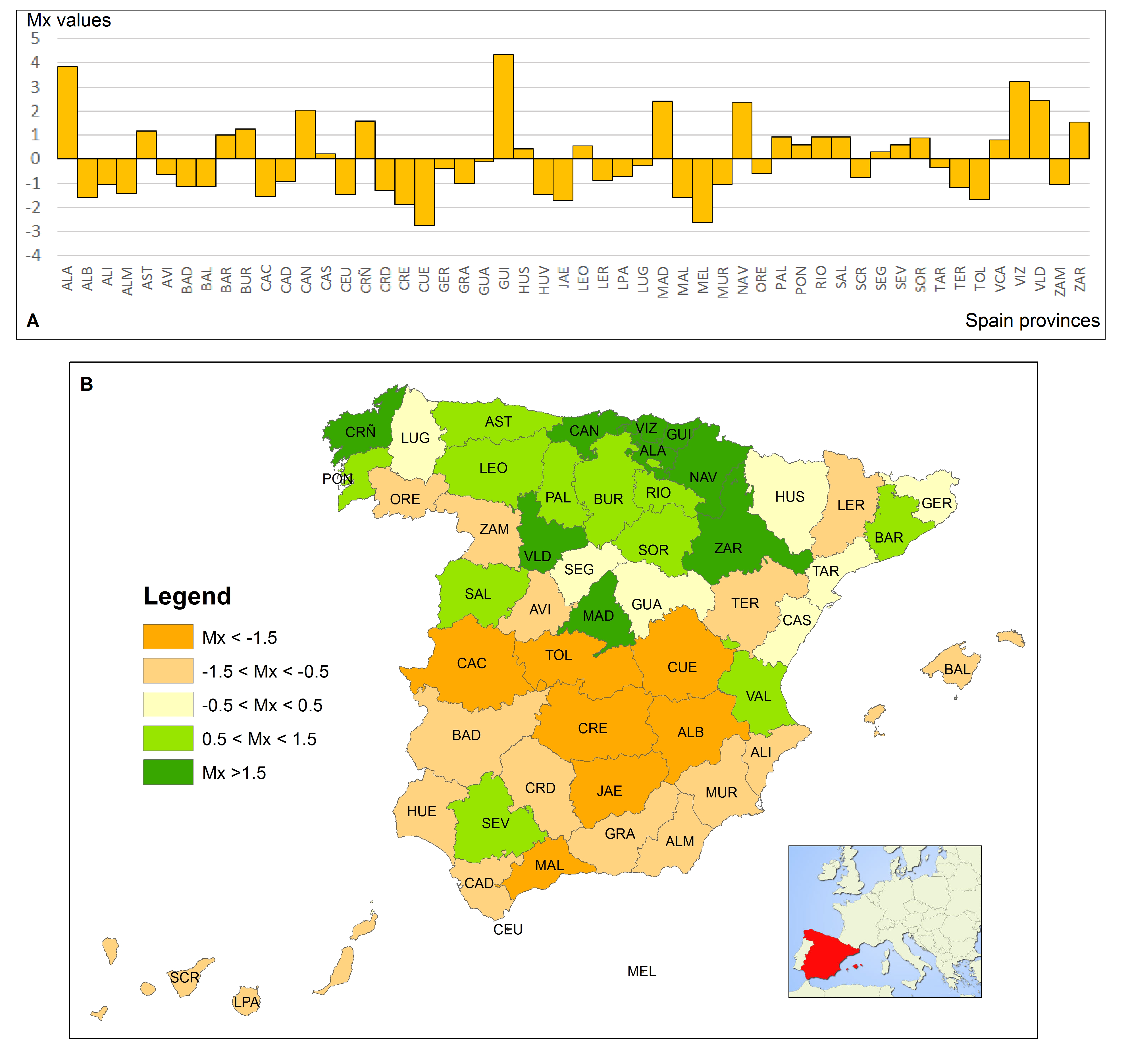} %
\end{center}
\caption{\emph{(A) Bar diagram with the values of $Mx$ obtained in each province of Spain. The optimal place to locate the private academy shall be where the value of Mx is the highest. (B) Geographic distribution of the $Mx$ values. Green provinces are the best places for our multiobjective problem. Design problem described in Equation \eqref{alb1}.}}
\end{figure}

\section{Conclusions}

\subsection{General conclusions} The generalized supporting vectors are proved to be useful in many applied disciplines. A large number of optimization problems in Physics, Engineering, Statistics, etc, can be reformulated in terms of generalized supporting vectors. The idea behind this reformulation is the {\em linearization} of the optimization problem. Once an optimization problem is {\em linearized}, then certain matrices come up and they are to be optimized by using matrix norms. Thus generalized supporting vectors arise. We then can exactly solve the following optimization problem: $$\left\{\begin{array}{l} \max \|A_1x\|_2^2+ \cdots + \|A_nx\|_2^2 \\ \|x\|_2=1\end{array}\right.$$ where $A_1,\dots, A_n$ are matrices. In fact, due to our Theorem \ref{gsv}, we can also solve exactly the infinite dimensional version of the previous optimization problem: $$\left\{\begin{array}{l} \max \sum_{i=1}^\infty\|T_i(x)\|^2 \\ \|x\|=1 \end{array}\right.$$ where the $T_i$’s are compact operators on Hilbert spaces such that $\sum_{i=1}^\infty \|T_i\|^2 <\infty$. We have shown by means of Subsection \ref{quantum} that the previous optimization problem has useful applications in the study of {\em effects} in Quantum Physics. A first application of supporting vectors was given in \cite{CSGPGRH} where a TMS coil was truly optimally designed. In that paper a three-component problem is stated but only the case of one component was solved. Here we solve the three-component case. Another application of supporting vectors is provided in this paper where an optimal location problem is solved using statistics multivariate analysis.

\subsection{Novelties of this work} In this subsection we enumerate the novelties provided by this work:
\begin{enumerate}
\item We provide an exact solution of an optimization problem, not a heuristic method for approaching it \cite{mente1,key:wassermann,key:romei}. Our exact solution is an analytical solution justified by a mathematical theorem whereas the heuristic methods presented in the literature are not mathematically proved to be convergent to the solution of the optimization problem.
\item For the first time in the literature of generalized supporting vectors, a MATLAB code is provided for computing them.
\item By means of the generalized supporting vectors we can exactly solve the three-component problem to obtain a truly optimal TMS coil, whereas until now the one-component problem was the only problem addressed \cite{CSGPGRH}. Our TMS coil is optimal because we solve the optimization problem exactly.
\item In virtue of our main result (Theorem \ref{gsv}), we are able to optimize a sequence of observable magnitudes by a pure state in a quantum mechanical system. In particular we focus on the probability density operator on a quantum mechanical system represented by an infinite dimensional separable complex Hilbert space. Hence, we find, by means of our Theorem \ref{gsv}, the state that maximizes the modulus of the probability density operator, which in fact coincides with the state of higher probability of the system.
\item We spotlight the relation of supporting vectors with Statistical Optimization, connecting multivariable analysis with multiobjective problem solutions.
\item This is an interdisciplinary work that comprises pure abstract nontrivial theorems with their proofs and programming codes with their results to directly apply them to real-life situations.
\end{enumerate}

\noindent {\bf Acknowledgements. } The authors would like to express their deepest gratitude towards the reviewers for their valuable suggestions and comments that helped improve the paper considerably.

\appendix

\section{Algorithms developed in this work}

In this section we show the algorithms written in MATLAB that allow to solve all the real problems that can be modeled with generalized supporting vectors (in particular, the problems presented in the previous two sections).

\subsection{General algorithm for generalized supporting vectors}

First we include the algorithm to compute the solution of the problem presented in \cite[Theorem 3.3]{CSGPMPSM}, that is
\begin{eqnarray*}
    \max_{\|x\|_2 = 1} \sum_{i=1}^k \| A_i x\|_2^2 &=& \lambda_{\max} \left(\sum_{i=1}^k A_i^T A_i\right) \\
    \arg\max_{\|x\|_2 = 1} \sum_{i=1}^k \| A_i x\|_2^2  &=& V\left(\lambda_{\max}\Big(\sum_{i=1}^k A_i^T A_i\Big)\right) \cap \E_{\ell_2^n}
\end{eqnarray*}
where $A_1,\dots,A_k$ are $m\times n$ real matrices. This algorithm beholds the particular case where $k=1$ and $A_1$ is an $m\times 2$ matrix whose column vectors have the same Euclidean norm (see Theorem \ref{teor3.5}).
\begin{lstlisting}
    function [lambda_max, x] = sol_1(M)
    %%%%%
    %%%%% INPUT:
    %%%%%
    %%%%% M = {A_1 A_2 ... A_k} a list with the matrices
    %%%%%
    %%%%%%%%%%%%%%%%%%%%%%%%%%%%%%%%%
    %%%%%
    %%%%% OUTPUTS:
    %%%%%
    %%%%% lambda_max - maximum eigenvalue
    %%%%% x - basis of unit  eigenvectors associated to lambda_max
    %%%%%%%%%%%%%%%%%%%%%%%%%%%%%%%%%
    %%%%
    k = length(M);
    [nrows ncols] = size(M{1});
    a1 = M{1}(:,1); %% First column of M{1}
    a2 = M{1}(:,2); %% Second column of M{1}
    if (k==1) & (ncols==2) & (abs(norm(a1)-norm(a2))<1e-12) 
        %%% In this particular case, a single matrix with two columns with 
        %%% the same norm is considered. A tolerance of 1e-12 is needed in 
        %%% order to compare these norms.
        %%% The maximum lambda_max and the supporting vectors are computed 
        %%% directly as in Theorem 3.5  
        lambda_max = norm(a1)^2 + abs(a1'*a2);
        if floor(a1'*a2)==0
            %%% The columns of this matrix form a basis of supporting vectors
            x = eye(2); 
        elseif a1'*a2>0
            %%% The columns of this matrix form a basis of supporting vectors
            x = [sqrt(2)/2 sqrt(2)/2; -sqrt(2)/2 -sqrt(2)/2]'; 
        elseif a1'*a2<0
            %%% The columns of this matrix form a basis of unit supporting vectors
            x = [-sqrt(2)/2 sqrt(2)/2; sqrt(2)/2 -sqrt(2)/2]';
        end
    else
        %%% This is the general case
        suma = zeros(ncols);
        for i=1:k
            suma = suma + M{i}'*M{i};
        end
        %%%
        [V,D] = eig(suma);           %%% Computing the eigensystem
        lambda_max = max(diag(D)); %%% This is the maximum eigenvalue
        N = size(D,1);
        %%% Now we find the indices where the elements of the diagonal of the
        %%% matrix D are equal (with a tolerance of 1e-12) to lambda_max  
        ind_lambda_max = find(abs(diag(D)-lambda_max*ones(N,1))<1e-12);
        x = V(:, ind_lambda_max); %%% %%% The columns of this matrix form a
        %%% basis of unit supporting vectors associated to the maximum eigenvalue 
    end
end
\end{lstlisting}

\subsection{Particular algorithm for the TMS coil}

Here we include the code to compute the solution of Problem \eqref{ffm}:
\begin{equation*} 
\left\{
\begin{array}{l}
\max  \left\| E_z \psi \right\| _2 \\
\max  \left\| E_y \psi \right\| _2 \\
\max  \left\| E_x \psi \right\| _2 \\
\min   \psi^T R \psi  
\end{array}
\right.
\end{equation*}

\begin{lstlisting}
function psi = sol_psi(Ex, Ey, Ez, R)

C=chol(R);                    % Cholesky's decomposition of matrix R = C' * C

A1=Ex*inv(C);                 % Change of variable to each matrix involved
A2=Ey*inv(C);                 
A3=Ez*inv(C);                 

[lambda, phi] = sol_1({A1 A2 A3}); % We apply the algorithm to obtain                                            % phi = C * psi
psi=inv(C)*phi(:,1);              % We undo all changes to compute                                     % the solution psi employing the                                     % first column in phi

end
\end{lstlisting}

\end{document}